\documentclass[reqno]{amsart}

\input xy
\xyoption{all}
\usepackage{epsfig}
\usepackage{color}
\usepackage{amsthm}
\usepackage{amssymb}
\usepackage{amsmath}
\usepackage{amscd}
\usepackage{amsopn}
\usepackage{url}
\usepackage{hyperref}\hypersetup{colorlinks}

\usepackage{stmaryrd}

\usepackage[T1]{fontenc}

\usepackage{color} %\textcolor{red}{Test}

\definecolor{darkred}{rgb}{1,0,0} %can change the intensity in [0,1]
\definecolor{darkgreen}{rgb}{0.0,0.8,0.0} 
\definecolor{darkblue}{rgb}{0,0,1}

\hypersetup{colorlinks,
linkcolor=darkblue,
filecolor=darkgreen,
urlcolor=darkred,
citecolor=darkgreen}

\numberwithin{equation}{section}
\newtheorem {Theorem}{Theorem}
\numberwithin{Theorem}{section}

\newtheorem {Lemma}[Theorem]{Lemma}

\newtheorem {Corollary}[Theorem]{Corollary}
\theoremstyle{definition}
\newtheorem{Definition}[Theorem]{Definition}
\theoremstyle{remark}
\newtheorem{Remark}[Theorem]{Remark}
\newtheorem{Example}[Theorem]{Example}

\expandafter\chardef\csname pre amssym.def at\endcsname=\the\catcode`\@
\catcode`\@=11
\def\undefine#1{\let#1\undefined}
\def\newsymbol#1#2#3#4#5{\let\next@\relax
 \ifnum#2=\@ne\let\next@\msafam@\else
 \ifnum#2=\tw@\let\next@\msbfam@\fi\fi
 \mathchardef#1="#3\next@#4#5}
\def\mathhexbox@#1#2#3{\relax
 \ifmmode\mathpalette{}{\m@th\mathchar"#1#2#3}%
 \else\leavevmode\hbox{$\m@th\mathchar"#1#2#3$}\fi}
\def\hexnumber@#1{\ifcase#1 0\or 1\or 2\or 3\or 4\or 5\or 6\or 7\or 8\or
 9\or A\or B\or C\or D\or E\or F\fi}

\font\teneufm=eufm10
\font\seveneufm=eufm7
\font\fiveeufm=eufm5
\newfam\eufmfam
\textfont\eufmfam=\teneufm
\scriptfont\eufmfam=\seveneufm
\scriptscriptfont\eufmfam=\fiveeufm

\catcode`\@=\csname pre amssym.def at\endcsname

\def    \eps    {\epsilon}

\newcommand{\A}{{\mathcal A}}
\newcommand{\tA}{\tilde{\mathcal A}}

\def    \R      {{\mathbb R}}

\def    \Z      {{\mathbb Z}}
\def    \N      {{\mathbb N}}
\def    \Q      {{\mathbb Q}}

\def    \T      {{\mathbb T}}
\def    \CP     {{\mathbb C}{\mathbb P}}

\def    \12    {{\frac{1}{2}}}

\def    \HC     {\operatorname{HC}}
\def    \CC     {\operatorname{CC}}

\def    \H      {\operatorname{H}}

\def  \MUCZ  {\operatorname{\mu_{\scriptscriptstyle{CZ}}}}

\begin{document}

%%%%%%%%%%%%%%%%%%%%%%%%%%%%%%
%   TEXT FORMATTING

\setlength{\smallskipamount}{6pt}
\setlength{\medskipamount}{10pt}
\setlength{\bigskipamount}{16pt}

%%%%%%%%%%%%%%%%%%%%%%%%%%

%%%%%%%%%%%%%%%%%%%%%%%%%%

%%%%%%%%%%%           BEGINNING OF  TEXT

%%%%%%%%%%%%%%%%%%%%%%%%%%

\title[Perfect Reeb flows]{Perfect Reeb flows and action-index relations}

\author[Ba\c sak G\"urel]{Ba\c sak Z. G\"urel}

\address{Department of Mathematics, University of Central Florida,
  Orlando, FL 32816, USA} \email{basak.gurel@ucf.edu}



\date{\today} 

\thanks{The work is partially supported by the NSF grant 
DMS-1414685.} %DMS-1207680

\begin{abstract} 
  We study non-degenerate Reeb flows arising from perfect contact
  forms, i.e., the forms with vanishing contact homology
  differential. In particular, we obtain upper bounds on the number of
  simple closed Reeb orbits for such forms on a variety of contact
  manifolds and certain action-index resonance relations for the
  standard contact sphere. Using these results, we reprove a theorem
  due to Bourgeois, Cieliebak and Ekholm characterizing perfect Reeb
  flows on the standard contact three-sphere as non-degenerate Reeb
  flows with exactly two simple closed orbits.
\end{abstract}

\maketitle

\tableofcontents
\section{Introduction and main results}
\label{sec:intro}
\subsection{Introduction}
\label{sec:introduction}
In this paper, we investigate non-degenerate Reeb flows arising from
contact forms with vanishing contact homology differential, referred
to as \emph{perfect} contact forms throughout. In particular, we
prove, under minor additional conditions, upper bounds on the number
of simple closed Reeb orbits for perfect forms on a broad class of
contact manifolds and certain action-index resonance relations for the
spheres $S^{2n-1}$ equipped with the standard contact structure.

To put this work in perspective, recall that a well-known conjecture
in Reeb dynamics asserts the existence of at least $n$ simple closed
Reeb orbits on the standard contact $S^{2n-1}$. (See, e.g.,
\cite{EH87,Lo,LZ,Wa} for the lower bounds on the number of such orbits
in the convex case and the discussion below for $S^3$.)  Furthermore,
hypothetically, the number of closed Reeb orbits on $S^{2n-1}$ is
either $n$ or infinite.  Moreover, drawing an analogy to several
results and conjectures concerning Hamiltonian flows, it is not
unreasonable to conjecture that there are infinitely many simple
closed Reeb orbits whenever there are more such orbits than necessary
(e.g., $n$ for $S^{2n-1}$).  We will refer to the latter conjecture,
which is the main motivation for the present work, as the
\emph{contact HZ-conjecture}, for, to the best of the author's
knowledge, the first written statement of the Hamiltonian
HZ-conjecture is in \cite{HZ}; see
\cite{BH,CKRTZ,Fr1,Fr2,FH,GG:hyperbolic,GHHM,Gu:hq,Gu:nc,Ke:JMD,LeC}
for some results in the Hamiltonian setting. All of these conjectures
are open in general, even under the non-degeneracy assumption.

In contrast with the Hamiltonian case, it is not clear how to express
the threshold number of simple closed orbits in the contact
HZ-conjecture in terms of the contact homology. However, one can
certainly interpret the condition that the Reeb flow has no
unnecessary closed Reeb orbits by requiring the contact differential
to vanish, i.e., the Reeb flow to be perfect. Note that, although
perfect flows are clearly exceptional (see Remark
\ref{rmk:GG-generic}), many Reeb flows of interest are in this class;
see Section \ref{sec:ex} for examples. For $S^3$ such flows were
studied in \cite{BCE}.

These conjectures lead to the question whether there is an upper bound
on the number of simple closed orbits of a perfect Reeb flow on, say,
the standard contact $S^{2n-1}$. In this work, we establish such upper
bounds of contact homological nature for the sphere and some other
contact manifolds under some additional assumptions; see Theorems
\ref{thm:even} and \ref{thm:all}. (In general, without extra
conditions, it is not even known if a perfect flow on $S^{2n-1}$,
$2n-1\geq 5$, must have finitely many simple closed orbits.) These
upper bounds, in particular, provide homological interpretation for
the upper bound $n$ in the HZ-conjecture for $S^{2n-1}$. To the best
our knowledge, these are the first results in this direction beyond
the three-dimensional case considered in \cite{BCE}.

Another aspect of our investigation of perfect Reeb flows is a
resonance relation asserting that all simple closed Reeb orbits $x$ of
a perfect Reeb flow on the standard contact $S^{2n-1}$ have the same
ratio $\Delta(x)/\A(x)$, where $\Delta(x)$ is the mean index and
$\A(x)$ is the action. This result, Theorem \ref{thm:res}, is a
contact analog of the action-index resonance relations proved in
\cite{CGG,GG:gaps} for Hamiltonian flows and is also related to the
results from \cite{EH87} and, at least on the conceptual level, to
\cite{HR}.

Using these results, we also give a short proof of a theorem from
\cite{BCE}, characterizing perfect Reeb flows on the standard contact
sphere $S^3$ as non-degenerate Reeb flows with exactly two simple
closed orbits; this is Theorem \ref{thm:BCE}.

It should be noted that much more is known about the number of simple
closed Reeb orbits in dimension three. This is partly due to the
availability of powerful, but strictly 3-dimensional, methods such as
finite energy foliations and the embedded contact homology and also
because the Conley-Zehnder index and local contact homology behave in
a much simpler way in this case.  For instance, the Weinstein
conjecture has been established in dimension three, \cite{Ta}, and,
moreover, it has been recently proved that a Reeb flow on any closed
contact 3-manifold has at least two closed orbits, \cite{CGH}; see
also \cite{GGo:index,GHHM,LL} for the case of the standard contact
$S^3$.  Furthermore, there exist either two or infinitely many closed
characteristics on any strictly convex hypersurface in $\R^4$ and,
under some additional assumptions,
% provided that all stable and unstable manifolds of the hyperbolic
% periodic orbits intersect transversally
for any non-degenerate Reeb flow on the standard contact $S^3$;
\cite{HWZ98,HWZ03}. (See also, e.g., \cite{Ho,HWZ96,WHL} and references
therein for some other relevant results.)

\subsection{Main results}
\label{sec:results}
The main object underlying our results is a perfect Reeb flow on a
contact manifold.
\begin{Definition}
\label{def:perfect}
  A non-degenerate contact form is said to be \emph{perfect} if the
  contact homology differential (for a chosen flavor of contact
  homology) vanishes. 
\end{Definition}
The flow of such a contact form will be called a perfect Reeb flow.
This notion, in general, depends on the kind of contact homology used,
but once this choice is made, perfectness is independent of the almost
complex structure and other auxiliary data. There are many examples
of perfect Reeb flows of interest; here, deferring further examples to
Section \ref{sec:ex}, we merely note that a non-degenerate Reeb flow
on an irrational ellipsoid is perfect. It is worth pointing out that
in all the examples known to us the Reeb flow is perfect even in a
stronger sense: all closed Reeb orbits in every free homotopy class
have the same Conley-Zehnder index modulo 2; let us call such Reeb
flows \emph{geometrically perfect}. Observe that a Reeb flow can be
geometrically perfect whether or not the contact homology is
defined. If, however, a flavor of contact homology is defined, a
geometrically perfect flow is also perfect in our sense. In what
follows ``perfectness'' for us will always be understood in the sense
of Definition \ref{def:perfect}.

Our first result is an analog of Theorem 1.1 from \cite{CGG} for Reeb
flows (see also \cite{GG:gaps}).  The Hamiltonian version of this
result asserts in particular that for a Hamiltonian diffeomorphism of
$\CP^n$ with finitely many periodic orbits, there are $n+1$ distinct
fixed points $x_0,\dotsc,x_n$ such that the so-called augmented
actions for all of these points are equal:
$\tA(x_0)=\cdots=\tA(x_n)$. (In general, the augmented action of a
contractible periodic orbit $x$ is defined by $\tA(x)=\A(x)- \lambda/2
\Delta(x)$, where $\A(x)$ and $\Delta(x)$ are, respectively, the
symplectic action and the mean index of $x$ and $\lambda$ is the
rationality constant of the ambient symplectic manifold.)  Similarly,
for perfect Reeb flows on the standard contact sphere $S^{2n-1}$ we
have

\begin{Theorem}
\label{thm:res}
Let $\alpha$ be a perfect non-degenerate contact form on
$(S^{2n-1},\xi_{std})$.  Then, for any two closed Reeb orbits $x$ and
$y$ of $\alpha$, we have
$$
\frac{\Delta(x)}{\A(x)}=\frac{\Delta(y)}{\A(y)}.
$$
\end{Theorem}

\begin{Remark}
  Notice that, in contrast with the Hamiltonian case, in Theorem
  \ref{thm:res} we make no assertion about the number of simple closed
  Reeb orbits.  At this stage we cannot rule out that the Reeb flow of
  $\alpha$ can have any number greater than or equal to two (including
  infinity) of simple closed orbits. (It is not hard to see that a
  non-degenerate contact form on the standard contact $S^{2n-1}$ must
  have at least two simple closed Reeb orbits; see Remark
  \ref{rmk:sphere} and also \cite{Kn}.)
\end{Remark}
 
As an immediate consequence of this theorem and a result from
\cite{BCE}, reproved below (see Theorem \ref{thm:BCE}), asserting that
a non-degenerate contact form on the standard contact $S^3$ is perfect
if and only if it has exactly two simple closed Reeb orbits, we obtain

\begin{Corollary} 
\label{cor:S3}
Let $\alpha$ be a non-degenerate contact form on $(S^3,\xi_{std})$
with exactly two simple closed Reeb orbits, say, $x$ and $y$.  Then
$\Delta(x)/\A(x) = \Delta(y)/\A(y)$.
\end{Corollary}

The next two theorems below hold for a more general class of contact
manifolds $(M^{2n-1},\xi)$ which we now specify. First, we require a
version of contact homology $\HC_*(M,\xi)$ to be defined for $M$. In
this paper we will utilize the linearized contact homology, although
cylindrical contact homology would equally well suit our
purposes. (See also Remark \ref{rmk:rigor}.)  To this end, assume that
$(M^{2n-1},\xi)$ is closed and strongly fillable by an exact
symplectically aspherical manifold, i.e., $M$ is the boundary of some
exact symplectic manifold $(W^{2n},\,\omega=d\tilde{\alpha})$ such
that $c_1(TW)=0$ and $\ker {\tilde{\alpha}}|_M = \xi$ with matching
orientations.
% The latter orientation condition simply means that the Liouville
% vector field $Y$ along $M$ determined by the requirement $i_Y \omega =
% {\tilde{\alpha}}|_M$ points outward.
Under these conditions, the linearized contact homology $\HC_*(M,\xi)$
is defined and independent of the contact form on $M$ supporting $\xi$
(see, e.g., \cite{Bo2,BO:inven,SFT}); we will use the notation
$\HC_*(\xi)$ whenever $M$ is clear from the context. Furthermore, we
will work with contact homology defined for the entire collection of
free homotopy classes of closed Reeb orbits, although everything below
would also go through if we restricted the homology to a smaller
collection which is closed under iterations; see Remark
\ref{rmk:homotopyclass}. (See Section \ref{sec:prelim} for a brief
discussion of the linearized contact homology.)
 
In addition, suppose that $\HC_*(\xi)$ satisfies the following
condition:
\begin{equation}
\label{eq:finiteness}
\dim \HC_* (\xi) < \infty  \text{ for all } * \leq 2n-4
\text{ and } \HC_* (\xi)=0 \text{ for } * \ll 0.
\end{equation}

This condition holds for a large variety of contact manifolds of
interest, but, of course, not for all contact manifolds. Next, set

\begin{equation}
\label{eq:b}
b =\limsup_{m\to\infty} \sum_{i=0}^{2n-2}\dim \HC_{m+i} (\xi).
\end{equation}

For instance, $b=n$ for the sphere $(S^{2n-1},\xi_{std})$. Also,
whenever $ \dim \HC_m (\xi) \leq b_0 $ with $b_0$ independent of the
degree $m$, we have $b\leq b_0 \cdot (2n-1)$ in general and $b\leq b_0
\cdot n$ if $\HC_m (\xi) \neq 0$ only in even degrees.
% $b\leq b_0 \cdot (n-1)$ if $\HC_* (M) \neq 0$ only in odd degrees.
Note that in \eqref{eq:b} we have also assumed that $\dim \HC_m (\xi)
< \infty$ for all large $m$. Otherwise, by definition, $b=\infty$ and
the results below hold trivially.  Finally, a simple closed orbit is
called \emph{even} if the number of real eigenvalues of the linearized
Poincar\'e return map in the interval $(-1,\, 0)$ is even; see Section
\ref{sec:prelim}. Otherwise, an orbit is said to be \emph{odd}.

Under the above hypotheses on $(M^{2n-1},\xi)$ with $b$ as in
\eqref{eq:b}, we have
\begin{Theorem}
\label{thm:even}
Let $\alpha$ be a perfect non-degenerate contact form on
$(M^{2n-1},\xi)$ with $r$ simple even orbits. Then $r \leq b$.
\end{Theorem}

\begin{Corollary} 
\label{cor:good}
A non-degenerate perfect Reeb flow on the standard contact sphere
$S^{2n-1}$ has at most $n$ even closed Reeb orbits.
\end{Corollary} 

Furthermore, considering \emph{all} simple closed orbits, we have

\begin{Theorem} 
\label{thm:all} 
Let $\alpha$ be a perfect non-degenerate contact form on
$(M^{2n-1},\xi)$. Assume that $\alpha$ has $r$ simple orbits (even or
odd) such that the reciprocals of their mean indices are linearly
independent over $\Q$. Then $r \leq b$.
\end{Theorem}

As an application of Theorem \ref{thm:even}, we reprove the following
result from \cite{BCE}.

\begin{Theorem}[\cite{BCE}]
\label{thm:BCE}
For a non-degenerate contact form $\alpha$ on $(S^3,\xi_{std})$, the
following are equivalent:
\begin{itemize}
\item[(i)] $\alpha$ is perfect.
\item[(ii)]$\alpha$ has exactly two simple closed Reeb orbits; these
  orbits are elliptic and their mean indices and actions are the same
  as those for an irrational ellipsoid. (These orbits are also
  unknotted and have linking number one.)
\item[(iii)] $\alpha$ has exactly two simple closed Reeb orbits.
\end{itemize}
\end{Theorem}

\begin{Remark}
  It is worth pointing out that it is unknown in general, even for the
  standard contact $S^{2n-1}$ when $2n-1 \geq 5$, whether a perfect
  Reeb flow necessarily has finitely many simple closed orbits.
\end{Remark}

\begin{Remark}
\label{rmk:rigor}
We should note that the paper relies on the theory of contact homology
which is still to be fully rigorously established; see
\cite{HWZ,HWZ2}. To circumvent the foundational difficulties, we could
have used the $S^1$-equivariant symplectic homology rather than the
linearized contact homology; see \cite{BO12}. However, the proofs in
that case would certainly be more involved and less transparent, and
hence we preferred to employ contact homology.
\end{Remark}

\begin{Remark}
\label{rmk:GG-generic}
Combining the argument from \cite{GG:generic} and Theorem
\ref{thm:all}, it is not hard to show that for any contact manifold
$(M,\xi)$ meeting the requirements of the theorem, non-perfect contact
forms form a residual subset in the set of all contact forms
supporting $\xi$, equipped with $C^{\infty}$-topology.
\end{Remark}

\begin{Remark}
\label{rmk:howcommon}
The linear independence condition from Theorem \ref{thm:all} is
satisfied for many perfect Reeb flows. Moreover, one can even
conjecture that the condition holds for a generic perfect Reeb flow;
for instance, this is the case for ellipsoids. Note, however, that the
resonance relation from \cite{GK09,Vi:res} asserts that the
reciprocals of mean indices and 1 are never linearly independent over
$\Q$ when the flow has finitely many simple closed orbits.
\end{Remark}

\subsection{Examples}
\label{sec:ex}
Most of the examples of perfect Reeb flows known to the author
originate from contact torus actions. In the most general setting,
whenever we have a contact torus action with only finitely many
one-dimensional orbits and such that there exists a one-parameter
subgroup in the torus with orbits transverse to the contact structure,
the action of such a generic one-parameter subgroup gives rise to a
perfect Reeb flow. (The closed Reeb orbits of this flow are exactly
the one-dimensional orbits of the torus action, and moreover, as is
easy to see, all these orbits are non-degenerate and elliptic.) Among
the contact manifolds this observation applies to are, for instance,
contact toric manifolds; see \cite{AM,Le}.  Moreover, one can further
perturb this Reeb flow, making it ergodic and creating non-trivial
dynamics by applying the results from \cite{Ka}.

Here, however, rather than exploring this construction in the most
general framework, we focus on more specific examples which we find
interesting and illuminating.

\begin{Example}[Prequantization circle bundles]
\label{ex:preq}

Let $(B,\omega)$ be a closed symplectic manifold equipped with a
Hamiltonian circle action with isolated fixed points and let $H$ be a
Hamiltonian (a moment map) generating the action. Assume in addition
that $B$ is integral, i.e., $\omega \in \H_2(B,\Z)/\mathrm{Tor}$. Let
$\pi \colon M \to B$ be the pre-quantization circle bundle over
$(B,\omega)$, i.e., $M$ is an $S^1$-bundle over $B$ with
$c_1=[\omega]$. In other words, $\pi^* \omega = d\alpha$, where
$\alpha$ is a connection form on $M$, and the Lie algebra of $S^1$ and
$\R$ are suitably identified. As is well-known, $\alpha$ is a contact
form. Denote by $a \in \R/\Z$ the action of $H$ on an orbit of the
action. (Here viewing $a$ as an element of $\R/\Z$ rather than $\R$
eliminates the ambiguity in the definition of the action arising from
the choice of a capping.)  It is easy to see that $a$ is well-defined,
i.e., independent of the orbit.

Let us now lift the $S^1$-action on $B$ to a contact $\R$-action on
$M$, preserving the form $\alpha$.  This lift is generated by the
vector field $H R + \hat{X}_H$. Here $R$ is the Reeb vector field of
$\alpha$, i.e., the vector field generating the principal bundle
$S^1$-action on $M$, and $\hat{X}_H$ is the horizontal lift of the
Hamiltonian vector field $X_H$. (We refer the reader to \cite[Appendix
A]{GGK} for a discussion of group actions and pre-quantization circle
bundles, including all the facts used in this example.)
Alternatively, when $H>0$, this vector field can be described as the
Reeb vector field of the form $H^{-1} \alpha$, where, abusing
notation, we write $H$ for the pull-back $H \circ \pi$.

One can show that when $a$ is irrational, i.e., $a \not\in \Q/\Z$, the
only simple periodic orbits of the resulting Reeb flow are the fibers
of $\pi$ over the fixed points of the $S^1$-action on $B$, i.e., the
critical points of $H$. Under our assumptions, there are finitely many
simple closed orbits and, in fact, the number of these orbits $r$ is
equal to the sum $\beta$ of the Betti numbers of $B$. Furthermore, all
such orbits are elliptic, since the critical points of $H$ are
elliptic, and hence the Reeb flow is geometrically perfect. The
requirements that $H>0$ and that $a \not\in \Q/\Z$ can always be
satisfied by adding a constant to $H$. Clearly, one can replace the
circle action by a torus action in this construction.

As a side remark, recall that when $a$ is rational, i.e., $a \in \Q/\Z
$, the lifted action is periodic, i.e., it is an $S^1$-action.
However, the minimal period $T$ of the lifted action need not be equal
to the period of the $S^1$-action on $B$. More specifically, when
$a=p/q$ with $(p,q)=1$ and the period of the $S^1$-action on $B$ is
one, we have $T=q$.
\end{Example}

Let us now briefly examine how Example \ref{ex:preq} fits in the
framework of Theorem \ref{thm:even}. In all examples of symplectic
manifolds admitting Hamiltonian circle actions with isolated fixed
points known to the author, $c_1(TB)\not=0$. Under this assumption,
$c_1(\xi)=0$, where $\xi=\ker \alpha$, only when $B$ is (strictly)
monotone or negative monotone. To fix a grading of $\HC_*(\xi)$, let
us choose a non-vanishing section $\mathfrak{s}$ of the square of the
complex determinant line bundle of $\xi$ (see Section
\ref{sec:prelim}). Then working with cylindrical contact homology and
using a suitably chosen sequence of perfect contact forms, we have
$$
\HC_*(\xi)=\bigoplus_{m=1}^\infty \H_{*-2+m\Delta}(B,\Q),
$$
where $\Delta$ is the mean index of the fiber of $\pi$ evaluated with
respect to $\mathfrak{s}$; see \cite{Bou02} and also \cite{KvK} and
references therein. (For example, when the fiber is contractible and
$B$ is monotone, $\Delta=2N$, where $N$ is the minimal Chern number of
$B$; \cite[p.\ 100]{Bou02}.)  Now, it is clear that $b=r=\beta$, the
sum of Betti numbers of $B$, when $\Delta> 2n-2$. When $\Delta\leq
2n-2$, we may have $b>r=\beta$. In the next example, we will see
explicitly that this indeed can be the case.

\begin{Example}[Katok--Ziller flows]
\label{ex:Ka}
The pre-quantization example described above includes, in particular,
the Katok-Ziller description of perfect geodesic flows for asymmetric
Finsler metrics on the spheres and some other manifolds. For the sake
of simplicity, let us focus on the sphere $S^n$.

In this case, $M = ST^*S^n$ is the unit cotangent bundle of $S^n$ and
$R$ is the spray of the round metric on $S^n$, and the geodesic flow
is clearly a circle action. The oriented geodesics are just the great
circles on $S^n$, and taking the quotient by the action, we obtain the
principle circle bundle $\pi\colon M\to B$ where
$B=\mathrm{Gr}^+(2,n+1)$, the Grassmannian of oriented 2-planes in
$\R^{n+1}$.

Consider now a linear $S^1$-action on $\R^{n+1}$. This action
naturally gives rise to a Hamiltonian circle action on $B$. When the
``weights'' of the $S^1$-action on $\R^{n+1}$ are relatively prime,
the fixed points of the resulting action on $B$ are isolated and,
choosing $H$ appropriately as in the previous example, we can lift
this action to a perfect Reeb flow on $M$.  One can show that this
Reeb flow is, in fact, the geodesic flow of an asymmetric Finsler
metric on $S^n$, described in \cite{Ka,Zi}. The number of simple
closed Reeb orbits in this case is $r=n+1$ when $n$ is odd and $r=n$
when $n$ is even, while $b=n+3$ when $n$ is odd and $b=n+2$ when $n$
is even. (See \cite{KvK} and references therein for the calculation of
contact homology in this case.) This example shows that the inequality
$r \leq b$ in Theorem \ref{thm:even} can be strict.

\end{Example}

\begin{Example}[GMSW--AM contact structures]
\label{ex:AM}
A sequence of geometrically perfect contact forms $\alpha_k$ on $S^2
\times S^3$ with exactly $r=4$ simple closed Reeb orbits is
constructed in \cite{AM,GMSW}. These contact forms have vanishing
first Chern class; their contact homology is concentrated in even
degrees and $\dim \HC_{2m}(\ker \alpha_k) = 2k +2$ when $2m>2$.  It
readily follows that $b=6(k+1)$, and hence the discrepancy between $r$
and $b$ in Theorem \ref{thm:even} can be arbitrarily large.
\end{Example}

\begin{Remark}
  A particular feature of the examples described above is that the
  Reeb flow generates a contact torus action. We emphasize that, in
  general, perfect Reeb flows need not have such a simple
  dynamics. For instance, the ergodic perturbations of perfect Finsler
  flows on $S^n$ constructed in \cite{Ka} are still perfect and have
  the same periodic orbits as the original flows, but exhibit highly
  non-trivial dynamics.
\end{Remark}

There are, however, examples of perfect Reeb flows of an entirely
different nature.

\begin{Example}[Hyperbolic geodesic and Reeb flows]
  The geodesic flow (or its $C^1$-small perturbation) on the unit
  cotangent bundle $M$ to a manifold $B$ with negative sectional
  curvature is geometrically perfect (see Section \ref{sec:results}).
  Indeed, there exists exactly one closed orbit in every free homotopy
  class of loops in $B$.  More generally, any hyperbolic Reeb flow is
  necessarily geometrically perfect; \cite{MP}. The reason is that the
  parity of the index of a hyperbolic orbit is determined by whether
  the orbit is even or odd, and an orbit is odd if and only if the
  restriction of the stable bundle to the orbit is not
  orientable. Now, it is clear that for a hyperbolic Reeb flow all
  orbits in a fixed free homotopy class have the same index mod 2.  We
  refer the reader to \cite{FoHa} for the constructions of ``exotic''
  hyperbolic Reeb flows.
\end{Example}

\begin{Remark}[Morse inequalities] The classical Morse inequalities
  give lower bounds for the number of critical points of a function in
  terms of the Betti numbers of a manifold, and the Morse inequalities
  turn into equalities when the Morse function is perfect. The
  non-degenerate Arnold conjecture, as we understand it now, is
  essentially a variant of the Morse inequalities for fixed points of
  Hamiltonian diffeomorphisms and Floer homology. One can wonder if
  there are analogs of the Morse inequalities for Reeb flows and
  contact homology. Obviously, providing a lower bound for all, not
  necessarily simple, closed Reeb orbits of a fixed index presents no
  difficulty.  However, once the count is restricted to simple orbits,
  the situation becomes much more subtle.  Asymptotic versions of
  contact Morse inequalities are established in \cite{GK09,
    GGo:index}; cf.\ \cite{HM} and also \cite[Proposition 3.4]{GG:nm}
  for a Hamiltonian variant of asymptotic Morse
  inequalities. (Strictly speaking, these results require the Reeb
  flow to have finitely many simple closed orbits, but one can expect
  the inequalities to hold, under some reasonable assumptions, without
  this requirement.)  But, even in the perfect case, there seems to be
  no obvious analog of the Morse (in)equalities. Theorems
  \ref{thm:even} and \ref{thm:all} give an inequality going in the
  opposite direction, and the examples of this section show that one
  can have arbitrarily large homology generated by very few simple
  orbits.
\end{Remark}

\subsection{Acknowledgements} The author is grateful to Alberto
Abbondandolo, Viktor Ginzburg, Nancy Hingston and Wolfgang Ziller for
useful discussions.

\section{Preliminaries}
\label{sec:prelim}

Assume that $(M^{2n-1},\xi)$ is a closed contact manifold strongly
fillable by an exact symplectically aspherical manifold, i.e.,
$M=\partial W$ for some exact symplectic manifold
$(W^{2n},\omega=d\tilde{\alpha})$ such that $c_1(TW)=0$ and $\ker
{\tilde{\alpha}}|_M = \xi$ with matching orientations. The latter
orientation condition simply means that the Liouville vector field $Y$
along $M$ determined by the requirement $i_Y \omega =
{\tilde{\alpha}}|_M$ points outward. As mentioned earlier, the
linearized contact homology $\HC_*(M,\xi)$ is then defined and
independent of the contact form on $M$ supporting $\xi$; see, e.g.,
\cite{Bo2,BO:inven,SFT}. We will write $\HC_*(\xi)$ whenever $M$ is
clear from the context. (Note that $\HC_*(\xi)$ also depends on the
filling, but for the sake of brevity this dependence is suppressed in
the notation.)

Let $\alpha$ be a non-degenerate contact form on $(M^{2n-1},\xi)$,
i.e., $\ker \alpha = \xi$ and all closed Reeb orbits of $\alpha$,
including the iterated ones, are non-degenerate.  A simple Reeb orbit
of $\alpha$ is called \emph{even} (\emph{odd}) if the number of real
eigenvalues of the linearized Poincar\'e return map in the interval
$(-1,\, 0)$ is even (odd). Moreover, an iterated orbit is said be
\emph{bad} if it is an even iteration of an odd periodic
orbit. Otherwise, an orbit is said to be \emph{good}. (For instance,
all simple orbits are good.) The contact homology $\HC_*(\xi)$ is the
homology of a complex $\CC_*(\alpha)$ generated over $\Q$ by good
closed Reeb orbits of $\alpha$. The contact homology complex is graded
by the Conley-Zehnder index of the Poincar\'e return map plus $n-3$,
i.e., the degree of a generator $x$ is $|x|=\MUCZ(x) + n-3$. The
differential, roughly speaking, counts the number of rigid punctured
holomorphic cylinders (with holomorphic cappings of the punctures in
$W$) in the symplectization; see, e.g., \cite{Bo2,BO:inven}. The
precise workings of the differential will be inessential for us since
we will mainly focus on perfect Reeb flows.

In the presence of non-contractible orbits, to have a well-defined
grading of the complex $\CC_*(\alpha)$ (or to evaluate the
Conley-Zehnder index of a periodic orbit) one needs an extra
structure. To this end, we fix a non-vanishing section $\mathfrak{s}$
of the square of the complex determinant line bundle of $TW$, which
exists since $c_1(TW)=0$; see, e.g., \cite{Es,GGo:index}.  Then, to
evaluate $\MUCZ(x)$ (or the mean index described below) of a periodic
orbit $x$, it suffices to choose a unitary trivialization of $x^* TW$
such that the square of its top complex wedge is
$\mathfrak{s}|_x$. Such a trivialization is unique up to homotopy, and
the Conley-Zehnder index (and the mean index) is well-defined.  The
grading of the contact homology complex $\CC_*(\alpha)$ and, in turn,
the homology $\HC_*(\xi)$ depends on the choice of the section
$\mathfrak{s}$, unless, of course, only the contractible orbits are
taken into account.

The mean index $\Delta(x)$ of a closed Reeb orbit $x$ (possibly
degenerate) measures the total rotation of certain eigenvalues on the
unit circle of the linearized Poincar\'e return map of the Reeb flow
at $x$; see \cite{Es,Lo,SZ}. When $x$ is non-degenerate,
\begin{equation}
\label{eq:Delta-MUCZ}
0\leq|\Delta(x)-\MUCZ(x)|<n-1.
\end{equation}
Furthermore, the mean index is homogeneous with respect to iteration:
\begin{equation}
\label{eq:index-hom}
\Delta(x ^k)=k \Delta(x).
\end{equation}

Finally, the action of a closed Reeb orbit $x$ is defined by $A(x) =
\int_x\alpha > 0$.

\begin{Remark}
\label{rmk:homotopyclass}
The contact homology splits into a direct sum of homology spaces over
free homotopy classes of loops in $W$.  Hence one has a contact
homology group defined for any collection of free homotopy classes.
(The grading is independent of the section $\mathfrak{s}$ for the
``contractible part'' of the contact homology, but, in general, it
does depend on $\mathfrak{s}$.) In what follows, we will work with the
entire contact homology; however, our results hold for the contact
homology defined for any collection of free homotopy classes closed
under iterations.
\end{Remark}

\section{Proofs}
\label{sec:proofs}
In this section we prove the main results of this paper.

\subsection{Proof of Theorem \ref{thm:res}}
First, observe that by our assumption that $\alpha$ is perfect,
$\CC_*(\alpha)=\HC_*(S^{2n-1}, \xi_{std})$, where
$$ 
\HC_*(S^{2n-1}, \xi_{std}) =
\begin{cases}
  \Q & \text{for $* \ge 2n-2$ and even,} \\
   0 & \text{otherwise};
\end{cases}
$$
see, e.g., \cite{Bo2}. Hence there is exactly one periodic orbit in
every even degree $\geq 2n-2$, generating the homology in that
degree. The closed Reeb orbits are, therefore, strictly ordered by the
degree and, in turn, by the Conley-Zehnder index.  Furthermore, these
orbits are also ordered strictly by the action: the action decreases
as one goes down from the generator of $\HC_*$ to the generator of
$\HC_{*-2}$. Indeed, denote by $x$ and $y$ the generators of $\HC_*$
and $\HC_{*-2}$, respectively. By the main theorem of \cite{BO:inven},
$x$ and $y$ are connected by a holomorphic curve in the
symplectization, and $\A(x) - \A(y) >0$ since the only holomorphic
curve having zero $\omega$-energy is the cylinder over a Reeb orbit;
see, e.g., \cite[Lemma 5.4]{BEHWZ}. The two orderings of the orbits
coincide and we have
$$
\MUCZ(x) > \MUCZ (y) \text{ iff } \A(x) > \A(y).
$$
This observation is the key point of the proof of Theorem
\ref{thm:res}.

Let $k>0$ be an integer such that $\A(y^k) > \A(x)$. We count the
iterations of $x$ ``between'' $x$ and $y^k$ in two ways: by the action
and by the index.  The result of counting by the action is
$$
\left \lfloor \frac{\A(y^k) - \A(x)}{\A(x)} \right \rfloor
= 
\left \lfloor k \frac{\A(y)}{\A(x)} \right \rfloor -1.
$$
On the other hand, using \eqref{eq:Delta-MUCZ}, the count by the index
is
$$
\left \lfloor \frac{\MUCZ(y^k) - \MUCZ (x)}{\MUCZ (x)} 
\right \rfloor
=
\left \lfloor \frac{\Delta(y^k) - \Delta (x)}{\Delta (x)} 
\right \rfloor + \mathrm{O}(1)
= 
\left \lfloor k \frac{\Delta (y)}{\Delta (x)} 
\right \rfloor -1 + \mathrm{O}(1),
$$
where O(1) stands for an error bounded (in absolute value) by a
constant independent of $k$. Equality of the two counts implies that
$$
\left \lfloor k \frac{\A(y)}{\A(x)} \right \rfloor
=
\left \lfloor k \frac{\Delta (y)}{\Delta (x)}
\right \rfloor + \mathrm{O}(1).
$$
Finally, dividing by $k$ and letting $k\to \infty$, we obtain
$$
\frac{\Delta(x)}{\A(x)}=\frac{\Delta(y)}{\A(y)}.
$$
$\hfill\square$
\begin{Remark}
  Theorem \ref{thm:res}, in particular, implies that the ordering of
  the closed Reeb orbits by the mean index is also strict and agrees
  with the orderings by the action and by the Conley-Zehnder index.
\end{Remark}

\subsection{Proofs of Theorems \ref{thm:even} and \ref{thm:all}}
\label{sec:pfs-good-all}
In this section we prove Theorems \ref{thm:even} and \ref{thm:all}.

\subsubsection{Proof of Theorem \ref{thm:even}}
\label{sec:pf-even}
Let $\alpha$ be a perfect non-degenerate contact form on
$(M^{2n-1},\xi)$ such that $\HC_*(\xi)$ satisfies
\eqref{eq:finiteness}. The key to the proof of Theorem \ref{thm:even}
is the following combinatorial lemma.

\begin{Lemma} 
\label{lem:key} 
Assume that $\alpha$ has a collection of $r$ geometrically distinct
periodic orbits $x_1,\dotsc, x_r$. Then for every sufficiently small
$\eps >0$, there exist positive integers $k_1, \dotsc, k_r \geq 1$
such that the iterated orbits $x_1^{k_1},\dotsc, x_r^{k_r}$ all have
degrees in an interval of length $2(n-1)+\eps$. \end{Lemma}

Theorem \ref{thm:even} immediately follows from Lemma \ref{lem:key}.
Indeed, note that all (even) orbits and their iterations contribute to
the homology since $\alpha$ is perfect. Furthermore, by the definition
of $b$ in \eqref{eq:b}, there are at most $b$ such orbits with indices
in an interval of length $2(n-1)+\eps$ comprising the $2n-1$ possible
degrees. Thus by Lemma \ref{lem:key} the total number of even orbits
cannot exceed $b$.

We will now prove Lemma \ref{lem:key}.

\subsubsection{Proof of Lemma \ref{lem:key}}
\label{sec:pf-key}
We begin by making a few observations. First, $\Delta_i := \Delta
(x_i) > 0$ for each $i$. Indeed, the fact that $\Delta_i \neq 0$
follows from our assumption that $\dim \HC_*(\xi) < \infty $ for $*
\leq 2n-4$. For, otherwise, by \eqref{eq:Delta-MUCZ}, all iterations
$x_i^k$ would have degrees in the bounded interval $(-2,\,2n-4)$,
resulting in infinite homology in some degree in this range. The
positivity, $\Delta_i > 0$, is a consequence of the assumptions that
$\HC_*(\xi) = 0 $ for $* \ll 0 $ and $\alpha$ is perfect, for all
iterations of an orbit with negative index contribute to the homology
in degrees $< 0 $.

Below we assume without loss of generality that $\Delta_1=\max_{i}
\Delta_i$ where $i=1,\dotsc,r$.  Set $k=(k_1,\dotsc,k_r)$ with $k_i
\in \Z$ and consider the $(r-1)$ linear forms $L_1, \dotsc, L_{r-1}$
in the variable $k \in \Z^r$ defined as \begin{align*} L_i (k) &
  =\Delta_1 k_1 - \Delta_{i+1} k_{i+1} \text{ for } i=1,\dotsc,r-1.
\end{align*}
In this setting, for any $\delta>0$ given, by the Minkowski
approximation theorem, \cite[Corollary 3, Appendix B]{Ca}, there
exists an integer vector $k \neq 0$ such that
\begin{equation}
\label{eq:Mink}
|L_i(k)| < \delta \text{ for all } i=1,\dotsc,r.
\end{equation}
Refining this assertion, we will show that if $\delta$ is sufficiently
small, there exists such $k \in \Z^r$ with components $k_i>0$ for all
$i=1,\dotsc,r$. To this end, choose $\delta <\min\{1,\Delta_1,\dotsc,
\Delta_r\}$ and let $k_j \neq 0$ be the first non-zero component of
$k$. Let us show that $k_j$ is necessarily $k_1$. Indeed, if $k_1=0$,
by \eqref{eq:Mink} and our choice of $\delta$, we have
$|L_i(k)|=|\Delta_{i+1}k_{i+1}| < \delta < |\Delta_{i+1}| $.
Therefore, $|k_{i+1}| < 1$ for all $i=1,\dotsc,r-1$ and, since
$k_{i+1} \in \Z$, we conclude that $k_{i+1}=0$ for all
$i=1,\dotsc,r-1$. As a result, $k=0$ since we also have $k_1=0$. This
contradiction shows that $k_1\neq 0$. In fact, $k_1$ can be chosen to
be positive. To see this, observe that $k$ satisfies \eqref{eq:Mink}
if and only if $-k$ does; thus, by replacing $k$ by $-k$ if necessary,
we can assume that $k_1 > 0$.
 
Now, given any $\eps>0$ satisfying $\eps <
\min\{1,\Delta_1,\dotsc,\Delta_r\}$, let $\delta=\eps/2$. Then the
above discussion implies that $|\Delta_1 k_1 - \Delta_i k_i| < \eps/2$
for all $k_i$. By our choice of $\eps$, since $k_1\geq 1$ and
$\Delta_1>0$, we have $\Delta_1 k_1 - \eps/2 >0$ and hence $\Delta_i
k_i> 0 $. This implies that $k_i > 0$ for all $i$ since
$\Delta_i>0$. Finally, note that we also have $|\Delta_i k_i -
\Delta_j k_j| < \eps$ for all $i,j \in \{1,\dotsc,r\}$, i.e., the
indices of the orbits $x_1^{k_1},\dotsc, x_r^{k_r}$ are at most $\eps$
apart from each other. Therefore, by \eqref {eq:Delta-MUCZ}, we
conclude that their degrees lie in an interval of length
$2(n-1)+\eps$.  $\hfill\square$

%%%% THEOREM ALL 

\subsubsection{Proof of Theorem \ref{thm:all}}
\label{sec:pf-all}
Note that, just as in the proof of Theorem \ref{thm:even}, we can
assume that mean indices $\Delta_i>0$ for all $i$, and we set
$\Delta_1=\max_i \Delta_i$. Clearly, it suffices to establish a
version of Lemma \ref{lem:key} applicable in the setting of Theorem
\ref{thm:all}, i.e., when odd orbits are also allowed.  Recall from
Section \ref{sec:prelim} that only the odd iterations of odd orbits
contribute to the complex; hence, for a small enough $\eps>0$, it
suffices to find odd positive integers $k_1,\dotsc,k_r$ such that
$|\Delta_1 k_1 - \Delta_i k_i| < \eps/2$.

To this end, set $k_i= 2 m_i + 1$ and write
$$
|\Delta_1 k_1 -\Delta_i k_i| = 2\Delta_i|m \theta_i - m_i - \alpha_i|,
$$
where $m=m_1$ and $\theta_i=
%\frac{\Delta_1}{\Delta_i}
\Delta_1/\Delta_i $ and $\alpha_i = 1/2 - \Delta_1/2\Delta_i $. Now,
the first step is to find $m \in \Z$ such that $ \| m \theta_i -
\alpha_i\| < \delta, $ where $\delta < \min_i \eps/4\Delta_i$ and the
norm $\| \cdot \|$ is the distance to the closest integer.  It is easy
to see that $1,\theta_2,\dotsc,\theta_r$ are linearly independent over
$\Q$, for $1/\Delta_i$'s are so by the hypotheses of Theorem
\ref{thm:all}. Therefore, by Kronecker's theorem, see, e.g.,
\cite{Ca}, the orbit $\{m\cdot(\theta_2,\dotsc,\theta_r) \mid
m\in\Z\}$ of $ (\theta_2,\dotsc,\theta_r)$ is dense in the torus
$\T^{r-1}$. In fact, both the positive and negative semi-orbits are
dense. Hence, we can find arbitrarily large $m \in \Z$ such that $m
\cdot (\theta_2,\dotsc,\theta_r)$ is sufficiently close to
$(\alpha_2,\dotsc,\alpha_r)$, i.e., $ \| m \theta_i - \alpha_i\| <
\delta $. This proves that there exist integers $m_i$ for
$i=1,\dotsc,r$ such that $|m \theta_i - m_i - \alpha_i| < \delta$.

To complete the proof, it suffices to show that $m_i$'s are positive
when $m$ is large. This is obvious since $\theta_i>0$ and we have $m_i
\geq m\theta_i - \alpha_i - \delta > 0$ for large $m$.
$\hfill\square$

%%%% THEOREM BCE

\subsection{Proof of Theorem \ref{thm:BCE}}
\label{sec:pf-bce}
To set the stage for the proof, we begin with a few observations. (We
emphasize that most of what follows holds only for $S^3$ and when all
closed Reeb orbits are non-degenerate.) Observe that there are two
types of closed Reeb orbits on $S^3$: elliptic and hyperbolic, and we
call a hyperbolic orbit negative (positive) hyperbolic if the
eigenvalues of the linearized Poincar\'e return map lie on the
negative (positive) real axis.  Clearly, elliptic and positive
hyperbolic orbits, denoted below by $x$, are even and negative
hyperbolic orbits, denoted below by $y$, are odd in the sense of
Section \ref{sec:prelim}.  Moreover, the Conley-Zehnder index is odd
for elliptic and negative hyperbolic orbits and even 
for a positive hyperbolic orbit.

Denote, as above, by $\Delta(z)$ the mean index of an orbit $z$. By
non-degeneracy, $\Delta(z) \not\in \Q$ if $z$ is elliptic.  When $z$
is negative (positive) hyperbolic, $\Delta(z) = \MUCZ(z)$ is an odd
(even) integer. Set
$$\sigma(z)=-(-1)^{\MUCZ(z)}.$$
So, $\sigma(z)$ is the Poincar\'e index of the orbit $z$, modulo the
factor $-1$. Note that $\sigma$ is $1$ for elliptic and negative
hyperbolic orbits and $-1$ for positive hyperbolic orbits.

% Set $$\sigma(z)=(-1)^{|z|}=-(-1)^n(-1)^{\MUCZ(z)}.$$ So, $\sigma(z)$
% is the topological (Poincar\'e) index of the orbit $z$, modulo the
% factor $-(-1)^n$.
% % or, more precisely, of the return map of $z$.
% On $S^3$, by the previous observations, $\sigma$ is $1$ for elliptic
% and negative hyperbolic orbits and $-1$ for positive hyperbolic
% orbits.

One of the key ingredients of the proof of Theorem \ref{thm:BCE} is a
resonance relation from \cite{GK09,Vi:res}, generalized to the
degenerate case in \cite{GGo:index,HM,WHL}. This relation holds in a
much greater generality, but we will only need it for non-degenerate
Reeb flows on the standard contact $S^3$. Let $\alpha$ be a
non-degenerate contact form on $(S^3, \xi_{std})$ with finitely many
simple periodic orbits. Then the resonance relation in the original
form as in \cite{Vi:res} is
\begin{equation}
\label{eq:res}
{\sum} 
\frac{\sigma(x_i)}{\Delta(x_i)}
+\frac{1}{2}{\sum}
\frac{\sigma(y_i)}{\Delta(y_i)}= 
\frac{1}{2}.
%\chi(S^3,\xi_{std}),
\end{equation}
Here the first sum runs over simple even (i.e., elliptic or positive
hyperbolic) orbits with positive mean index and the second sum runs
over simple odd (i.e., negative hyperbolic) orbits with positive mean
index.

Finally, recall that the degree of a closed Reeb orbit $x$ in the
contact homology complex for $(S^3, \xi_{std})$ is $|x|=\MUCZ(x)-1$,
and
$$ 
\HC_*(S^3, \xi_{std}) =
\begin{cases}
  \Q & \text{for $* \geq 2$ and even,} \\
   0 & \text{otherwise}.
\end{cases}
$$

We now turn to the proof of Theorem \ref{thm:BCE}.

\subsubsection*{(i) $\Rightarrow$ (ii).} Assume that $\alpha$ is a
non-degenerate perfect contact form on $S^3$ supporting the standard
contact structure. 

Observe first that $\alpha$ has no positive hyperbolic orbits, for
$\HC_*(S^3,\xi_{std})$ is non-vanishing only in even degrees (i.e.,
when the Conley-Zehnder index is odd) and $\alpha$ is perfect.
Secondly, by Theorem \ref{thm:even}, there are at most two simple
elliptic orbits.

As above, denote by $x$ the even simple orbits and by $y$ the odd
simple orbits of $\alpha$. So, $x$'s are elliptic and $y$'s are
negative hyperbolic. Recall from the proof of Theorem \ref{thm:even}
that all mean indices are positive.  We claim that there is at most
one negative hyperbolic orbit of $\alpha$. To prove this claim, assume
the contrary and let $y_1$ and $y_2$ be two distinct negative
hyperbolic orbits with $\MUCZ(y_i) = \Delta(y_i) = 2 m_i + 1$ for
$i=1,2$ and $m_i\geq 1$. Then, using \ref{eq:index-hom}, we have
$$
\MUCZ(y_1^{2 m_2 +1}) = (2 m_1 +1)(2 m_2 +1) = \MUCZ(y_2^{2 m_1 +1}).
$$
In other words, $\alpha$ has two distinct orbits contributing to the
homology in degree $m=(2 m_1 +1)(2 m_2 +1)-1$, which contradicts the
fact that $\HC_m(S^3,\xi_{std})=\Q$.

To recap, we have established that $\alpha$ has finitely many simple
periodic orbits, and among these at most two are elliptic and at most
one is negative hyperbolic.  Moreover, when $\alpha$ is perfect, the
relation \eqref{eq:res} is
\begin{equation}
\label{eq:res1}
\sum_{x_i} \frac{1}{\Delta(x_i)}
+\frac{1}{2}\sum_{y_i}\frac{1}{\Delta(y_i)}= \frac{1}{2},
\end{equation}
where $x_i$'s are elliptic and $y_i$'s are negative hyperbolic.

Let us now show that $\alpha$ has at least two simple periodic
orbits. Indeed, assume that there is only one orbit, $z$.  Clearly,
$z$ cannot be elliptic by \eqref{eq:res1}, for elliptic orbits have
irrational mean index.  Hence $z$ must be negative hyperbolic with
$\MUCZ(z)=3$ and should generate the entire homology. However, it is
easy to see using \eqref{eq:Delta-MUCZ} that no odd iteration of $z$
can generate the homology $\HC_4(S^3,\xi_{std}) = \Q$.  Therefore,
$\alpha$ must have at least two simple closed Reeb orbits.  It is
worth pointing out that this assertion also follows from the more
general and highly non-trivial results of \cite{CGH} or from
\cite{GHHM,GGo:index}; see also Remark \ref{rmk:sphere}.

At this point there are three possibilities; the composition of simple
closed orbits of $\alpha$ could be two elliptic, one elliptic and one
negative hyperbolic or two elliptic and one negative hyperbolic. We
will next rule out the latter two alternatives by arguing that
presence of two different types of orbits is impossible.  To this end,
assume that $\alpha$ has (at least) one elliptic orbit, $x$, and one
negative hyperbolic orbit, $y$. Recall that $\Delta (x) \not \in \Q$
and $\Delta (y) = 2m+1 \in \N$. (Though this will follow from the
discussion below, notice that the second alternative is actually
impossible by \eqref{eq:res1}.)  Set $\theta = \Delta(x)/2\Delta(y)
\not \in \Q$.  Then, for any $\delta>0$, by Kronecker's theorem, see,
e.g., \cite{HaWr}, there exist $k, l \in \N$ such that $|k \theta - l
- 1/2| < \delta$. Therefore, for any $\eps>0$, choosing $\delta <
\eps/2\Delta(y)$, we obtain
\begin{equation*}
|k \Delta(x)- (2l+1) \Delta(y)| = 
|\Delta(x^k)- \Delta(y^{2l+1})|
< \eps.
\end{equation*}
Of course, this is equivalent to $ |\Delta(x^k)-\MUCZ(y^{2l+1})| <
\eps$ and, by \eqref{eq:Delta-MUCZ}, implies that
\begin{equation}
\label{eq:ellhyp}
|\MUCZ(x^k)-\MUCZ(y^{2l+1})| < 1+\eps.
\end{equation}
Now, recall that both indices in \eqref{eq:ellhyp} are odd.  Since the
distance between two distinct odd integers is at least two, we have
$\MUCZ(x^k)=\MUCZ(y^{2l+1})$, provided that $\eps<1$.  But then the
contact homology in degree $(2l+1)(2m+1)-1$ would be two dimensional,
which is impossible.

By the above discussion, $\alpha$ has exactly two simple closed Reeb
orbits and both of these orbits are elliptic.  The indices of these
orbits satisfy \eqref{eq:res1}, and all pairs of positive irrational
numbers satisfying this relation clearly occur as the mean indices for
irrational ellipsoids (cf.\ \cite{Lo}).  Moreover, Corollary
\ref{cor:S3} implies that the actions are also the same as those for a
suitable irrational ellipsoid.  The assertion that the orbits are
unknotted and have linking number one follows from the machinery of
\cite{HWZ98} as in \cite{BCE}.
% we are not claiming that the flows would equivalent, which is
% probably wrong. just claiming that only the same pairs as for an
% irrational ellipsoid can occur. Fayad-Katok examples, \cite{FK}.

\subsubsection*{(iii) $\Rightarrow$ (i).} Assume that $\alpha$ has
exactly two simple closed Reeb orbits. To conclude that $\alpha$ is
perfect, it suffices to show that both orbits have the same parity (of
Conley-Zehnder indices).  Thus we need to prove that $\alpha$ cannot
have simultaneously elliptic or negative hyperbolic orbits and
positive hyperbolic orbits.  The former case is impossible by
\eqref{eq:res}.

To see that the presence of one negative and one positive hyperbolic
orbit is impossible, assume the contrary and let $x$ be the positive
hyperbolic orbit and $y$ be the negative hyperbolic orbit.  Then
$\MUCZ(x)=\Delta(x)=2k$ and $\MUCZ(y)=\Delta(y)=2m+1$ for some
non-zero integers $k$ and $m$. Now, since $\HC_*(S^3,\xi_{std})$ is
non-vanishing only in positive even degrees, both mean indices must be
positive, and hence $k>0$ and $m\geq0$. Then the relation
\eqref{eq:res} reads
\begin{equation*}
\label{eq:res2}
-\frac{1}{2k}+\frac{1}{2(2m+1)}
= \frac{1}{2}.
\end{equation*}
In other words, $1/(2m+1) - 1/k = 1$, which is clearly impossible.
This proves that both of the orbits have the same parity, and hence
$\alpha$ must be perfect.

Theorem \ref{thm:BCE} follows from the implications (i) $\Rightarrow$
(ii) and (iii) $\Rightarrow$ (i) along with the obvious implication
(ii) $\Rightarrow$ (iii). $\hfill\square$

\begin{Remark}
\label{rmk:sphere}
A non-degenerate contact form $\alpha$ supporting the standard contact
structure on $S^{2n-1}$ must have at least two distinct simple closed
Reeb orbits. Indeed, assume that the Reeb flow of $\alpha$ has just
one simple closed (non-degenerate) orbit $x$.  Then we necessarily
have $\MUCZ(x)=n+1$ and $\Delta(x)=2$. Now, by \eqref{eq:Delta-MUCZ},
the orbit $x$ must be degenerate which contradicts the the background
non-degeneracy assumption.  (In fact, $x$ is a symplectically
degenerate maximum, \cite[Proposition 3]{GHHM}, and it also follows
from the results of that paper that in this case the Reeb flow must
have infinitely many periodic orbits.)
\end{Remark}

\end{document}